\input{style/arxiv-ba.cfg}
\documentclass[ba,linksfromyear,preprint]{imsart}
\makeatletter
   \@ifpackageloaded{natbib}{}{\usepackage{natbib}}
\makeatother
\usepackage{amsthm,amsmath}

\pubyear{2015}
\volume{10}
\issue{2}
\firstpage{505}
\lastpage{509}
\doi{10.1214/15-BA942B}

\begin{document}

\begin{frontmatter}
\title{Comment on Article by Dawid and Musio\thanksref{T1,T2}}
\runtitle{Comment on Article by Dawid and Musio}

\relateddois{T1}{Main article DOI: \relateddoi[ms=BA942]{Related item:}{10.1214/15-BA942}.}
\thankstext{T2}{This work is supported by the National Science
Foundation under award numbers DMS-1310294, SES-1024709, and
SES-1424481.}

\begin{aug}
\author[a]{\fnms{Christopher M.} \snm{Hans}\ead[label=e1]{hans@stat.osu.edu}}
\and
\author[b]{\fnms{Mario} \snm{Peruggia}\corref{}\ead[label=e2]{peruggia@stat.osu.edu}}

\runauthor{C. M. Hans and M. Peruggia}

\address[a]{Department of Statistics,
The Ohio State University,
Columbus, Ohio, U.S.A.,
\printead{e1}}

\address[b]{Department of Statistics,
The Ohio State University,
Columbus, Ohio, U.S.A.,\\
\printead{e2}}

\end{aug}

\end{frontmatter}


Dawid and Musio present interesting results on how to affect
model comparison using proper scoring rules, focusing chiefly on
Bayesian model comparison.
Among the reasons stated to justify the proposed
approach we note:
\begin{enumerate}
\item
The insensitivity of the procedure to a renormalization of the prior
distribution,
\item
The flexibility and/or robustness of the method when implemented
using a prequential score.
\end{enumerate}

The focus of the article is on the derivation of consistency results
for the proper scoring rule methods based both on their implementation
through a multivariate score and a prequential score.  There are very many such results in the article, but the gist of the argument is
that some form of proper scoring rule method can produce a consistent
procedure even in cases when the standard Bayesian approach fails to
do so or when it fails altogether, as is the case when improper
priors are used and Bayes factors cannot be calculated.

Consistent model selection is unquestionably a desirable property as
is the formulation of a coherent, universal framework for statistical
inference.  The Bayesian approach using {\em proper\/} priors
accomplishes the latter.  The proposed proper scoring rule methods
mend the complications that arise when the Bayesian approach is used
with improper priors.  However, the beauty of the coherent Bayesian
inferential framework is lost when model comparison is no longer based
on the likelihood score.  As in all compromises, something is gained
at the expense of losing something else, or, as some would
say, there is no free lunch!

Then, for those situations in which the Bayesian approach is not
broken, two questions arise naturally:
\begin{enumerate}
\item
When does a proper scoring rule model comparison produce a different
answer than a log-score model comparison?
\item
For those situations in which the answers are different, can an
argument be made for preferring the proper scoring rule method?
\end{enumerate}
This suggests juxtaposing the proposed method to model comparison
methods that compare directly the (log-) likelihoods for the various
models.\vadjust{\goodbreak}

Focusing on the technically simpler situations, such as that of the
univariate Gaussian process of Section~6.1, may be helpful to develop
some deeper intuition.  The addenda in the cumulative prequential
delta log-scores of~(7) in the article are given by
\[
S_{L,i}(x_i,Q_i) - S_{L,i}(x_i,P_i) =  \frac{1}{2}
\left[\log \sigma^2_{Q_i} - \log \sigma^2_{P_i} +
(x_i - \mu_{Q_i})^2/\sigma^2_{Q_i}
-(x_i - \mu_{P_i})^2/\sigma^2_{P_i}\right],
\]
and the addenda in the cumulative prequential delta Hyv\"arinen scores
are given by
\[
S_{H,i}(x_i,Q_i) - S_{H,i}(x_i,P_i) =
2/\sigma^2_{P_i} - 2/\sigma^2_{Q_i} +
(x_i - \mu_{Q_i})^2/\sigma^4_{Q_i} - (x_i -
  \mu_{P_i})^2/\sigma^4_{P_i},
\]
where $(\mu_{P_i}, \sigma^2_{P_i})$
and $(\mu_{Q_i}, \sigma^2_{Q_i})$ are the conditional means and variances
of $x_i$ given $\mathbf{x}^{i-1}$ (all the observations preceding
$x_i$), under models $P$ and $Q$, respectively.

Note that, for the case of a covariance stationary Gaussian process,
$\sigma^2_{P_i}$ and $\sigma^2_{Q_i}$ are constant in $i$.  As a
consequence, the cumulative prequential delta scores based on the
Hyv\"arinen rule and the log-score are perfectly linearly related
whenever $\sigma^2_{P_i} = \sigma^2_{Q_i} = \tau^2$.  As an example,
this is the case for two iid sequences with possibly different means
and equal variances and for two zero-mean, AR(1) sequences with
possibly different autoregressive parameters and equal innovation
variances.  The delta scores are also perfectly linearly related if the two
covariance stationary Gaussian process have equal conditional means
$\mu_{P_i} = \mu_{Q_i}$ and possibly different conditional variances
$\sigma^2_{P_i} = \tau^2_P$ and $\sigma^2_{Q_i} = \tau^2_Q$.  As an
example, this is the case for two iid sequences with equal means and
possibly different variances and for two zero-mean, AR(1) sequences
with equal autoregressive parameters and possibly different innovation
variances.

For Gaussian processes with non-stationary covariance structure, the
prequential delta scores based on the Hyv\"arinen rule and on the
log-score may not be perfectly linearly related.  Is it then possible to
characterize with necessary and sufficient conditions the Gaussian
processes for which the two delta scores are perfectly linearly
related?  When the delta scores are not perfectly linearly related,
how do they differ both in a finite-sample and an asymptotic sense?

Regardless of whether the delta scores are or are not perfectly
linearly related, there remains the question of how model comparison
decisions based on the two scores differ.  To address this issue, we
look at comparisons between two models and conform to the
recommendation made by the authors in Section~4, which is to select
the model with the lower prequential score.  When using the log-score,
this corresponds to using the Bayes decision rule under 0--1 loss and
assuming equal prior probabilities for the two models.  When using the
Hyv\"arinen score, there does not appear to be any principled way to
justify the use of the zero cut-off for the difference in prequential
scores, although, if the delta log-score and the delta Hyv\"arinen
scores are perfectly linearly related, such a cut-off for the delta
Hyv\"arinen score can be readily made to correspond to infinitely many
Bayes rules under generalized 0--1 loss for an appropriate choice of
prior model probabilities.

An inspection of the expressions for $S_{L,i}(x_i,Q_i) -
S_{L,i}(x_i,P_i)$ and $S_{H,i}(x_i,Q_i) - S_{H,i}(x_i,P_i)$ reveals
that the squared departures of the observations from their conditional
means are normalized by the conditional variance in the log-score and
by the {\em squared\/} conditional variance in the Hyv\"arinen score.
Beside the unnatural fact that the normalized terms are no longer
unitless, this suggests that the delta Hyv\"arinen score may be more
sensitive than the delta log-score to the presence of outlying
observations when the alternative model has larger variance than the
model generating the data.

\begin{figure}[t]
\centering
   \includegraphics{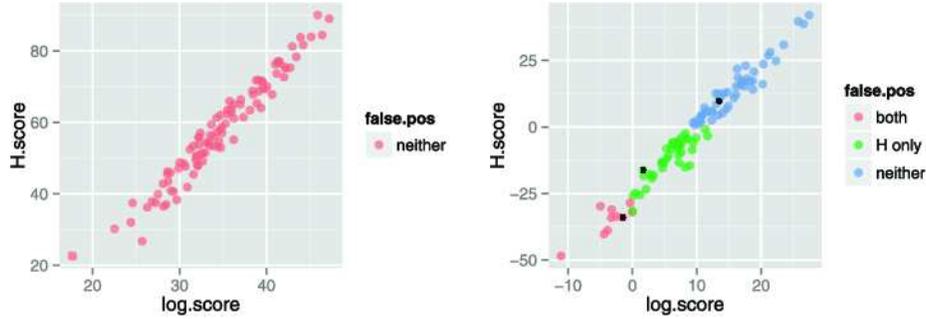}
   \caption{Cumulative prequential delta Hyv\"arinen scores vs.\ delta
     log-scores for 100 simulated data sets.  False positive
     identifications of $Q$ as the data generating model are
     highlighted by color.  The three points in the right panel
     plotted with a black center correspond to the three time series
     in Figure~\ref{fig:sample_ts}.}
  \label{fig:false_pos}
\end{figure}

This point is illustrated in Figure~\ref{fig:false_pos}, which is
based on 100 simulated data sets of size 101 from a zero-mean Gaussian
AR(1) process $P$ with autoregressive parameter $\phi$ equal to 0.5
and innovation variance equal to 1.  The alternative model, $Q$, is
taken to be a zero-mean Gaussian AR(1) process with autoregressive
parameter $\phi$ equal to 0.1 and innovation variance equal to 4.  The
prequential delta log-scores and delta Hyv\"arinen scores are built
based on the conditional distributions of observations 2 through 101.
These distributions are Gaussian with mean equal to $\phi$ times the
preceding observation and variance equal to the innovation variance.

For each simulated data set, we calculate the cumulative prequential
delta scores
\[
\Delta_L^{101}(\mathbf{x}^{101};P,Q)= \sum_{i=2}^{101} (S_{L,i}(x_i,Q_i) -
S_{L,i}(x_i,P_i))
\]
and
\[
\Delta_H^{101}(\mathbf{x}^{101};P,Q)=
\sum_{i=2}^{101}(S_{H,i}(x_i,Q_i) - S_{H,i}(x_i,P_i)).
\]
Correct identification of the data generating model under score $\ast$
corresponds to
\[
\Delta_{\ast}^{101}(\mathbf{x}^{101};P,Q) > 0.
\]

The left panel of Figure~\ref{fig:false_pos} shows that both the delta
log-score and the delta Hyv\"arinen score correctly identify model $P$
as the data generating model in all 100 simulations.  The right panel
corresponds to the same simulated data sets with the exception that,
in each data set, the 50th observation in the sequence of 101 is
contaminated by adding 7 to it, making the observation an additive
outlier.  The figure shows that the delta Hyv\"arinen score is much
more sensitive to the presence of the additive outlier.  In 10 out of
100 cases both methods incorrectly select model $Q$, in 50 cases they
both correctly select model $P$, but there are 40 cases in which only
the method based on the delta Hyv\"arinen score incorrectly selects
model $Q$.

\begin{figure}[t]
\centering
   \includegraphics{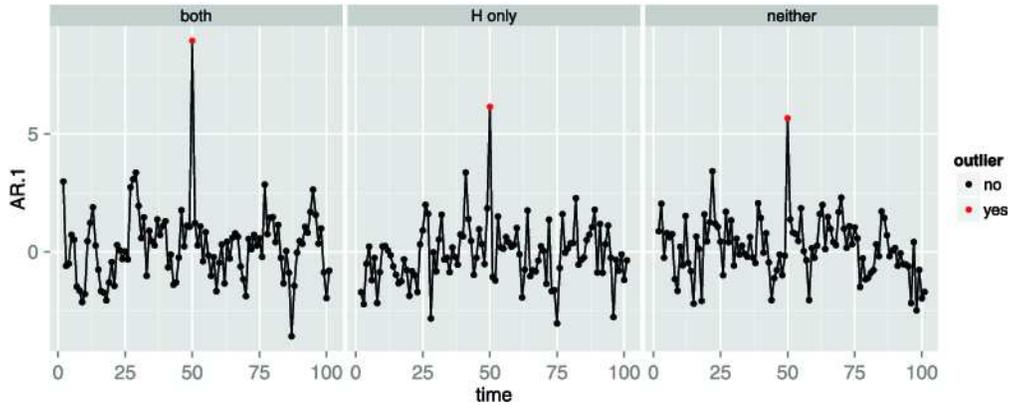}
   \caption{Three sample contaminated times series (with the outlier
     depicted in red) and their misclassification status according to
     the cumulative prequential delta Hyv\"arinen scores and delta log-scores.
     These three series correspond to the points plotted with
     a black center in the right panel of Figure~\ref{fig:false_pos}.}
  \label{fig:sample_ts}
\end{figure}

Figure~\ref{fig:sample_ts} displays three sample contaminated times
series (with the outlier depicted in red) that were analyzed in the
simulations.  The first series is misclassified by both delta scores,
the second is misclassified by the delta Hyv\"arinen score only, and
the third is correctly classified by both delta scores.  These three
series correspond to the points plotted with a black center in the
right panel of Figure~\ref{fig:false_pos}.

The normalization by the square of the variance (or an estimate of the
variance) appears throughout the article (cf. (34), (43),
and (70)) leading one to suspect that in all these situations the
Hyv\"arinen score may similarly be impacted by the presence of
outliers.  While the dependencies in the data may have played some
role in our simulation, we are convinced that the variance
normalization is the main issue.  In fact, we were able to simulate
examples with similar features after setting $\phi$ equal to zero in
both processes, thus eschewing the effect of serial correlations.
Such a choice, however, causes the delta scores to be perfectly
linearly related and makes the figures harder to decipher due to
overplotting, which is why we presented the simulation based on
correlated data instead.

Our simulation, following the set-up of Section~6.1, is based solely
on a comparison of the likelihoods for the two models.  However, it is
reasonable to conjecture that the sensitivity of the Hyv\"arinen score
to the presence of outliers will be injected, via the likelihood, also
when Bayesian model comparisons are carried out, irrespective of the
type of prior distribution specified for the model parameters (proper,
improper, subjective, or objective, as the case might be).  Related
questions are as follows.  Is it possible to modify the prequential
Hyv\"arinen score so as to alleviate its sensitivity to the presence
of outliers?  How does the method behave in the face of other model
violations?  Are other model comparison methods based on different
proper scoring rules not as sensitive to the presence of outliers?

In summary, the authors have proposed an interesting method for
performing Bayes\-ian model selection when improper priors are used for
within-model parameters by replacing the log marginal likelihood with
a proper scoring rule. The method avoids the machinations associated
with several of the alternative approaches that the authors mention
toward the end of Section~2 at the expense of moving even farther away
from the formal Bayesian paradigm.  The authors justify their approach
in part by proving consistency for model selection in certain
settings.

While the paper provides a framework for approaching the problem,
important choices still must be made in order to implement the
strategy, both with proper and improper priors.  We have seen that
these choices can have a substantial impact on the finite-sample
properties of the methods.  Our investigation was limited to the
Hyv\"arinen score, as this is the score most thoroughly discussed in
the paper.  The authors note that they ``confined attention to the
most basic homogeneous rule, the Hyv\"arinen score'' for simplicity,
but that ``there are no clear theoretical grounds for preferring one
[homogeneous scoring rule] over another.''  In light of our
investigation above, we wonder whether some theoretical progress might
be made by identifying a limited set of properties that might be of
interest (e.g., scale invariance, robustness to model violations,
etc.) and identifying classes of scoring rules and variants of the
prequential score that perform appropriately with respect to one or
all of these considerations. We believe that further research in this
direction would give the framework a stronger theoretical footing and
provide guidance to practitioners who wish to use the methods.

\end{document}